\documentclass[12pt]{amsart}
\usepackage{amssymb,amsfonts,amsmath,amsopn,amstext,amscd,latexsym,xy,color}

\input xy
\xyoption{all}

\theoremstyle{remark}

\newtheorem{para}{\bf}[section]

\theoremstyle{definition}

\theoremstyle{plain}

\newtheorem{thm}[para]{Theorem}

\newtheorem{prop}[para]{Proposition}

\newenvironment{numequation}{\addtocounter{para}{1}
\begin{equation}}{\end{equation}}

\newcommand{\vpi}{{\varpi}}

\newcommand{\vep}{{\varepsilon}}

\newcommand{\bbF}{{\mathbb F}}

\newcommand{\bbX}{{\mathbb X}}

\newcommand{\bbZ}{{\mathbb Z}}

\renewcommand{\frm}{{\mathfrak m}}

\newcommand{\fro}{{\mathfrak o}}
\newcommand{\frp}{{\mathfrak p}}

\newcommand{\cM}{{\mathcal M}}

\newcommand{\ra}{\rightarrow}
\newcommand{\lra}{\longrightarrow}

\newcommand{\sub}{\subset}

\newcommand{\midc}{{\; | \;}}

\newcommand{\fronr}{{\hat{\fro}^{nr}}}

\newcommand{\End}{{\rm End}}

\newcommand{\Aut}{{\rm Aut}}

\newcommand{\Spec}{{\rm Spec}}
\newcommand{\Spf}{{\rm Spf}}

\begin{document}

\begin{center}\Large

{\bf Galois actions on torsion points of universal one-dimensional
formal modules}

\bigskip

\large

Matthias Strauch

\bigskip
\normalsize

{\it
Department of Pure Mathematics and Mathematical Statistics\\
Centre for Mathematical Sciences, University of Cambridge\\
Wilberforce Road, Cambridge, CB3 0WB, United Kingdom\\
M.Strauch@dpmms.cam.ac.uk }

\vskip25pt

\end{center}

\small {\bf Abstract.} Let $F$ be a local non-Archimedean field
with ring of integers $\fro$ and uniformizer $\vpi$. Let $\bbX$ be
a one-dimensional formal $\fro$-module of $F$-height $n$ over the
algebraic closure $\bbF$ of the residue field of $\fro$. By the
work of Drinfeld, the universal deformation $X$ of $\bbX$ is a
formal group over a power series ring $R_0$ in $n-1$ variables
over the completion of the maximal unramified extension $\fronr$
of $\fro$. For $h \in \{0,\ldots,n-1\}$ let $U_h \sub \Spec(R_0)$
be the locus where the connected part of the associated
$\vpi$-divisible module $X[\vpi^\infty]$ has height $h$. Using the
theory of Drinfeld level structures we show that the
representation of $\pi_1(U_h)$ on the Tate module of the \'etale
quotient is surjective.

\normalsize

\tableofcontents

\section{Preliminaries}

Let $F$ be a local field non-Archimedean field with ring of
integers $\fro$ and uniformizer $\vpi$. Let $q$ be the cardinality
of the residue field of $\fro$, and denote this field by $\bbF_q$.
Fix a one-dimensional formal $\fro$-module $\bbX$ of $F$-height
$n$ over the algebraic closure $\bbF$ of the residue field of
$\fro$. By the work of Drinfeld, cf. \cite{D}, the universal
deformation $X$ of $\bbX$ is a formal $\fro$-module over a power
series ring $R_0$ in $n-1$ variables over $\fronr$, the completion
of the maximal unramified extension of $\fro$. We denote by

$$[ \cdot ]_X: \fro \lra \End_{R_0}(X)$$

\medskip

the action of $\fro$ on the universal deformation.

\bigskip

One may choose coordinates $u_1,\ldots,u_{n-1}$ of $R_0$ and a
coordinate $T$ on $X$ such that the multiplication by $\vpi$ on
$X$ is given by a power series $[\vpi]_X(T) \in R_0[[T]]$ with the
property that for all $i = 0,\ldots,n$:

\begin{numequation}\label{mult by vpi}
[\vpi]_X(T) \equiv u_iT^{q^i} \mbox{ mod } (u_0,\ldots,u_{i-1})
\,,  \mbox{ deg }(q^i + 1) \,,
\end{numequation}

\medskip

where we have put $u_0 = \vpi$ and $u_n = 1$ (cf. \cite{Ha}, sec.
21.5, \cite{HG}, Prop. 5.7).

\bigskip

The subset of $\Spec(R_0)$ where the $\vpi$-divisible
$\fro$-module $X[\vpi^\infty]$ is \'etale is precisely $\Spec(R_0)
- V(\vpi)$. Then the subset $U \sub \Spec(R_0/(\vpi))$ where the
\'etale quotient of $X[\vpi^\infty]$ has height $n-1$ is
$\Spec(R_0/(\vpi)) - V(u_1)$.\footnote{This is the locus where
$X[\vpi^\infty]$ is 'ordinary', where 'ordinary' means in this
context that the connected component of $X[\vpi^\infty]$ over $U$
is a Lubin-Tate formal $\vpi$-module of height one.} Denote by
$\kappa = \bbF((u_1,\ldots,u_{n-1}))$ the field of fractions of
$R_0/(\vpi)$ and put $\eta = \Spec(\kappa)$. Let $\kappa^a$ be an
algebraic closure of $\kappa$ and put $\bar{\eta} =
\Spec(\kappa^a)$. We consider the Tate module of the \'etale
quotient of $X[\vpi^\infty]$ over $\eta$:

$$T_1 := T_\vpi((X[\vpi^\infty]_\eta)^{\acute{e}t}) =
\lim_{\stackrel{\longleftarrow}{m}} X[\vpi^m]^{\acute{e}t}_\eta
(\kappa^a)$$

\medskip

which is a free $\fro$-module of rank $n-1$. The absolute Galois
group $\pi_1(\eta,\bar{\eta})$ of $\kappa$ acts on this
$\fro$-module $\fro$-linearly, and this action factors through
$\pi_1(U,\bar{\eta})$:

\begin{numequation}\label{conjecture}
\pi_1(\eta,\bar{\eta}) \lra  \pi_1(U,\bar{\eta}) \lra
\Aut_\fro(T_1)^\times \simeq GL_{n-1}(\fro) \,.
\end{numequation}

\medskip

In \cite{T}, Conj. 1.4, Y. Tian conjectures that this
representation is surjective. We give here a simple proof of the
surjectivity, using only results about Drinfeld level structures.
Tian's conjecture is more general than our result since it applies
to elementary $p$-divisble groups of any slope in $(0,1)$, not
only to one-dimensional groups of slope $\frac{1}{n}$ as is the
case considered here. In \cite{T} a proof of the surjectivity is
given the case of an elementary $p$-divisible group of slope
$\frac{1}{3}$. In \cite{B2} P. Boyer proves the irreducibility of
Igusa varieties in the case of certain 'simple' Shimura varieties
which had been studied previously in \cite{HT}. The irreducibility
of the 'Igusa varieties of the first kind' implies the
surjectivity of the Galois representation considered here (in
fact, it is equivalent to the surjectivitiy). Boyer derives this
result from a detailed study of the cohomology of these Shimura
varieties, cf. \cite{B1}. As already stated above, the method used
here employs only results about rings representing deformation
functors with level structures.

\medskip

Shortly after a first draft of this paper was written, Eike Lau
used the results presented here to prove the surjectivity of the
monodromy representation for certain Newton strata in the
universal deformation space of an arbitrary connected
Barsotti-Tate group over an algebraically closed field of
characteristic $p$. Lau's results cover in particular all cases of
Tian's conjecture, cf. \cite{L}.

\bigskip

\section{Statement of the result}

In this paper we prove a statement which is slightly more general
than the surjectivity of (\ref{conjecture}). For $h \in
\{0,\ldots,n-1\}$ put

$$R_{h,0} = R_0/(\vpi,u_1,\ldots,u_{h-1}) \,,$$

\medskip

with the convention that $u_0 = \vpi$, and so $R_{0,0} = R_0$.
Then the closed reduced subscheme of $\Spec(R_0)$ where the height
of the connected component of $X[\vpi^\infty]$ is at least $h$ is
equal to $\Spec(R_{h,0})$, and the open part of $\Spec(R_{h,0})$
where the height of the connected component is equal to $h$ is

$$U_h := \Spec(R_{h,0}) - V(u_h) \,.$$

\medskip

Hence the scheme $U$ considered above is equal to $U_1$. Let
$\kappa_h$ be the field of fractions of $R_{h,0}$ and put $\eta_h
= \Spec(\kappa_h)$. Let $\kappa_h^a$ be an algebraic closure of
$\kappa_h$ and put $\bar{\eta}_h = \Spec(\kappa_h^a)$. Fix a
positive integer $m$. Denote by

$$T_{h,m} := X[\vpi^m]^{\acute{e}t}_{\eta_h}(\kappa_h^a)$$

\medskip

the module of $\kappa_h^a$-valued points of the group scheme
$X[\vpi^m]^{\acute{e}t}_{\eta_h}$ over $\eta_h$. It is a free
$\fro/(\vpi^m)$-module of rank $n-h$. The absolute Galois group
$\pi_1(\eta_h,\bar{\eta}_h)$ of $\kappa_h$ acts $\fro$-linearly on
$T_{h,m}$, and we denote this representation by $\sigma_{h,m}$:

$$\sigma_{h,m}: \pi_1(\eta_h,\bar{\eta}_h) \lra
\Aut_\fro(T_{h,m}) \simeq GL_{n-h}(\fro/(\vpi^m)) \,.$$

\medskip

$\sigma_{h,m}$ clearly factors as $\pi_1(\eta_h,\bar{\eta}_h) \ra
\pi_1(U_h,\bar{\eta}_h) \ra \Aut_\fro(T_{h,m})$.

\medskip

\begin{thm}\label{result} For all $h \in \{0,\ldots,n-1\}$ and $m>0$ the
homomorphism $\sigma_{h,m}$ is surjective. In particular, denoting
by

$$T_h := \lim_{\stackrel{\longleftarrow}{m}}
X[\vpi^m]^{\acute{e}t}_{\eta_h}(\kappa_h^a)$$

\medskip

the Tate module of the \'etale quotient of
$X[\vpi^\infty]_{\eta_h}$, the resulting representation

$$\sigma_h: \pi_1(\eta_h,\bar{\eta}_h) \lra \Aut_\fro(T_h) \simeq
GL_{n-h}(\fro)$$

\medskip

is surjective.
\end{thm}

\bigskip

\section{Level structures}

Fix an integer $m>0$. Let $\cM_m = \Spf(R_m)$ be the formal scheme
over $\fronr$ which parameterizes deformations of $\bbX$ equipped
with a level-$m$-Drinfeld structure. It is proven in \cite{D},
Prop. 4.3,  that the ring $R_m$ is a regular local ring and a
finite flat $R_0$-algebra. Let

$$\phi := \phi^{univ}_m: (\vpi^{-m}\fro/\fro)^n \lra \frm_{R_m}$$

\medskip

be the universal level-$m$-structure. Here $\frm_{R_m}$ denotes
the maximal ideal of $R_m$ which is given the structure of an
$\fro$-module after having fixed a coordinate $T$ on the universal
deformation $X$. That $\phi$ is a level-$m$-structure means that
there is an invertible power series $\vep_m(T) \in R_m[[T]]$ such
that there is an equality

\begin{numequation}\label{eq1} \prod_{a \in
(\vpi^{-1}\fro/\fro)^n} (T-\phi(a)) = \vep_m(T) \cdot [\vpi](T)
\end{numequation}

\medskip

of formal power series over $R_m$, cf. \cite{D}, definition before
Prop. 4.3. For ease of notation put

$$A_m = (\vpi^{-m}\fro/\fro)^n$$

\medskip

and let $e_i = (0,\ldots,0,\vpi^{-m},0,\ldots,0) \in A_m$, be the
$i^{\mbox{\small{th}}}$ standard generator of $A_m$ as an
$\fro/(\vpi^m)$-module. By \cite{D}, Prop. 4.3, the sequence
$\phi(e_1),\ldots,\phi(e_n)$ is a regular system of parameters for
$R_m$, and the analogous fact is true for any basis of $A_m$. The
group

$$G_m := \Aut_\fro(A_m) \simeq GL_n(\fro/(\vpi^m))$$

\medskip

acts on $R_m$ via its action on the universal Drinfeld basis.
Concretely: $g \in \Aut_\fro(A_m)$ maps $\phi(a)$ to $\phi(g(a))$
for any $a \in A_m$. Let $\kappa_m = Frac(R_m)$ be the field of
fractions of $R_m$. Then $\kappa_m$ is a Galois extension of
$\kappa_0 = Frac(R_0)$ and the action of $G_m$ on $\kappa_m$
defines an isomorphism

\begin{numequation}\label{fund iso}
G_m \stackrel{\simeq}{\lra} Gal(\kappa_m/\kappa_0) \,.
\end{numequation}

(For details see \cite{St}, Thm. 2.1.2.) This means that the field
$\kappa_m$ generated by the $\vpi^m$-torsion points of $X$ is a
Galois extension of $\kappa_0$ with Galois group isomorphic to
$G_m$. (This was proven for the case $\fro = \bbZ_p$ by a
different method in \cite{RZ}.)

\bigskip

\section{Proof of the surjectivity}

Fix $m>0$ as above and $h \in \{0,\ldots,n-1\}$. We consider the
prime ideal $\frp_h = (\vpi,u_1,\ldots,u_{h-1})$ of $R_0$, and
determine the prime ideals of $R_m$ lying over $\frp_h$, as well
as the corresponding decomposition and inertia groups. Let

$$\frp_{h,m} = (\phi_1,\ldots,\phi_h) \sub R_m$$

\medskip

be the prime ideal generated by the regular sequence
$\phi_1,\ldots,\phi_h$, and put $R_{h,m} = R_m/\frp_{h,m}$,
$\kappa_{h,m} = Frac(R_{h,m})$. Denote by

$$A_{h,m} = \fro/(\vpi^m)e_1 \oplus \cdots \oplus
\fro/(\vpi^m)e_h \,.$$

\medskip

the submodule of $A_m$ generated by the first $h$ standard
generators. Define $P_{h,m} \sub G_m$ to be the stabilizer of the
submodule $A_{h,m} \sub A_m$ and let $Q_{h,m} \sub P_{h,m}$ be the
subgroup which acts trivially on the quotient $A_m/A_{h,m}$.

\medskip

\begin{prop}\label{decomposition} i) The ideal $\frp_{h,m} \sub
R_m$ lies over $\frp_h$, and the finite field extension
$\kappa_{h,m}/\kappa_h$ is normal.

ii) Via the isomorphism \ref{fund iso} the group $P_{h,m}$ is the
decomposition group of $\frp_{h,m}$.

iii) The inertia subgroup of $P_{h,m}$ is equal to $Q_{h,m}$, and
the canonical map

$$P_{h,m}/Q_{h,m} \lra \Aut(\kappa_{h,m}/\kappa_h)$$

\medskip

is a bijection. In particular, if $\kappa_{h,m}^s$ denotes the
maximal separable extension of $\kappa_h$ contained in
$\kappa_{h,m}$, then $Gal(\kappa_{h,m}^s/\kappa_h)$ is isomorphic
to $GL_{n-h}(\fro/(\vpi^m))$.
\end{prop}

\medskip

{\it Proof.} When $h=0$ the statements of this proposition are
equivalent to the assertion that $\kappa_m/\kappa_0$ is a Galois
extension with group $G_m$, cf. \ref{fund iso}. Since this is
known (\cite{St}, Thm. 2.1.2) we assume from now on $1 \le h \le
n-1$.

\medskip

i) We consider first the reduction of the left side of \ref{eq1}
modulo $\frp_{h,m}$ and see that as polynomials over
$R_m/\frp_{h,m}$ we have

\begin{numequation}\label{reduction}
\prod_{a \in \vpi^{m-1}A_m} (T-\phi(a)) \mbox{ mod } \frp_{h,m} \;
= \; \prod_{a' \in \vpi^{m-1}A'_{h,m}}(T-\phi(a'))^{q^h} \,,
\end{numequation}

where $A'_{h,m} = \fro/(\vpi^m)e_{h+1} \oplus \cdots \oplus
\fro/(\vpi^m)e_n$. Using \ref{mult by vpi} we see that modulo
$\frp_{h,m}$ the coefficients $\vpi,u_1,\ldots,u_{h-1}$ of
$[\vpi]_X(T)$ vanish, and so $\frp_h$ is necessarily contained in
$\frp_{h,m} \cap R_0$. Since both prime ideals have height $h$ and
$R_m$ is integral over $R_0$, we get that $\frp_{h,m} \cap R_0 =
\frp_h$. In particular, $u_h$ is not in $\frp_{h,m}$. Now consider
the coefficient of $T^{q^h}$ in \ref{reduction}. This coefficient
is equal to

$$\pm \prod_{a' \in \vpi^{m-1}A'_{h,m} - \{0\}} \phi(a') \,.$$

\medskip

Comparing this with the reduction of \ref{mult by vpi} modulo
$\frp_h$ shows that $\phi(a') \notin \frp_{h,m}$ for every
non-zero $a' \in \vpi^{m-1}A'_{h,m}$. A fortiori: $\phi(a') \notin
\frp_{h,m}$ for every non-zero $a' \in A'_{h,m}$. Since
$\phi(a+a') = \phi(a) +_X \phi(a') = \phi(a) + \phi(a') +
\phi(a)\left( \cdots \right)$, we deduce for later use that

\begin{numequation}\label{non-vanishing}
\mbox{For all } a' \in A_m - A_{h,m}: \; \phi(a') \notin
\frp_{h,m}
\end{numequation}

By \cite{Bou}, Ch. V, \S 2.2, Thm. 2, the extension
$\kappa_{h,m}/\kappa_h$ is normal.

\medskip

ii) $P_{h,m}$ is clearly contained in the decomposition group of
$\frp_{h,m}$. If $g \in G_m$ is not in $P_{h,m}$, then there is $i
\in \{1,\ldots,h\}$ such that $g(e_i)$ is not in $A_{h,m}$. Hence
$g$, considered as Galois automorphism, maps $\phi_i$ to an
element $\phi(g(e_i))$ which is not in $\frp_{h,m}$, by what we
have seen above, cf. \ref{non-vanishing}. This shows that
$P_{h,m}$ is the decomposition group of $\frp_{h,m}$.

\medskip

iii) The ring $R_{h,m}$ is a regular local ring with regular
system of parameters given by the images $\bar{\phi}_j$ of
$\phi_j$ for $j = h+1,\ldots,n$. The action of an element $g \in
P_{h,m}$ on $R_{h,m}$ sends $\bar{\phi}_j$ to $\phi(g(e_j)) \mbox{
mod } \frp_{h,m}$. In order for $g$ to induce the identity on
$R_{h,m}$ it is necessary and sufficient that

$$\phi(g(e_j)) \mbox{ mod } \frp_{h,m} \; = \;
\phi(e_j) \mbox{ mod } \frp_{h,m} \mbox{ for all } j =
h+1,\ldots,n \,.$$

\medskip

When we consider both sides of this equation as torsion points of
the formal group $X$ over $R_{h,m}$ we find that their difference,
as an element of the maximal ideal of $R_{h,m}$ (with the group
law induced by $X$), is equal to

$$\phi(g(e_j)) -_X \phi(e_j) \mbox{ mod } \frp_{h,m}
\; = \; \phi(g(e_j) -e_j) \mbox{ mod } \frp_{h,m} \,.$$

\medskip

If this difference is zero in $R_{h,m}$ then, by
\ref{non-vanishing}, we have $g(e_j)-e_j \in A_{h,m}$ for all $j =
h+1,\ldots,n$. This means that $g$ is in $Q_{h,m}$. The last
assertion is \cite{Bou}, Ch. V, \S 2.3, Prop. 6. \hfill $\Box$

\bigskip

{\bf Proof of Theorem \ref{result}.} We let $A'_{h,m} =
\fro/(\vpi^m)e_{h+1} \oplus \cdots \oplus \fro/(\vpi^m)e_n$ as
above. It follows from \ref{non-vanishing} that the induced map

$$\bar{\phi}: A'_{h,m} \lra R_{h,m} \sub \kappa_{h,m}$$

\medskip

is injective. Since $X[\vpi^m]^{\acute{e}t}_{\eta_h}$ has rank
$q^{m(n-h)}$ we see that $\kappa_{h,m}$ is obtained by adjoining
all $\vpi^m$-torsion points of $X_{\eta_h}$ to $\kappa_h$. Denote
by $\kappa_{h,m}^s$ the maximal separable extension of $\kappa_h$
contained in $\kappa_{h,m}$. Then the canonical maps

$$X[\vpi^m]_{\eta_h}^{\acute{e}t}(\kappa_{h,m}^s) \ra
X[\vpi^m]_{\eta_h}^{\acute{e}t}(\kappa_{h,m}) \leftarrow
X[\vpi^m]_{\eta_h}(\kappa_{h,m})$$

\medskip

are bijections. The connected component $X[\vpi^m]^\circ_{\eta_h}$
is infinitesimal of degree $q^{mh}$, and so the torsion points of
$X[\vpi^m]^{\acute{e}t}_{\eta_h}$ are exactly the elements
$\bar{\phi}(a')^{q^{mh}} \in \kappa_{h,m}$ which lie actually in
the maximal separable sub-extension $\kappa_{h,m}^s \sub
\kappa_{h,m}$. The action of the group

$$\Aut(\kappa_{h,m}/\kappa_h) = P_{h,m}/Q_{h,m}
\simeq \Aut_\fro(A'_{h,m}) \simeq GL_{n-h}(\fro/(\vpi^m))$$

\medskip

on the set of torsion points

$$\{ \bar{\phi}(a') \midc a' \in A'_{h,m} \}$$

\medskip

of $X[\vpi^m]_{\eta_h}$ induces an action on the set of torsion
points

$$\{ \bar{\phi}(a')^{q^{mh}} \midc a' \in A'_{h,m} \}$$

\medskip

of $X[\vpi^m]^{\acute{e}t}_{\eta_h}$. The latter action is again
equal to the full automorphism group $\Aut_\fro(T_{h,m})$. Hence
the map $Gal(\kappa_{h,m}^s/\kappa_h) \ra \Aut_\fro(T_{h,m})$ is
an isomorphism. \hfill $\Box$

\end{document}